\documentclass[final,12pt]{elsarticle}
\usepackage{xcolor,soul,cancel}
\usepackage[color]{showkeys}
\usepackage{amssymb,amsmath,amsthm,mathtools}
\usepackage{graphicx,float,enumitem,hyperref}
\newlist{alphalist}{enumerate}{1} 
\setlist[alphalist]{label=(\alph*)} 
\journal{arXiv}
\definecolor{refkey}{rgb}{0,1,1}
\definecolor{labelkey}{rgb}{1,0,0}

\numberwithin{equation}{section}

\usepackage[a4paper,margin=27mm]{geometry}
\usepackage{amsmath,amssymb,amsthm,mathtools}
\usepackage{microtype}
\usepackage{enumitem}

\newtheorem{thm}{Theorem}[section]
\newtheorem{lem}[thm]{Lemma}

\newtheorem{prop}[thm]{Proposition}

\newtheorem{defi}[thm]{Definition}

\theoremstyle{remark}

\newcommand{\R}{\mathbb R}
\newcommand{\E}{\mathbb E}

\newcommand{\Cov}{\operatorname{Cov}}

\newcommand{\Gini}{\mathcal D}

\theoremstyle{definition}

\begin{document}
\begin{frontmatter}
\title{A Sharp Perimeter Inequality for Mean Distances on Convex Curves}
\author[eric]{Eric Shen}
\ead{erick.2013@yandex.ru}
\address[eric]{Moscow State University, Moscow, 119991, Russia \\
Institute for Numerical Mathematics, Russian Academy of Sciences, Russia}
\date{August 2026}

\begin{abstract}

For a planar convex body $K$, let $\theta(K)$ denote the
expected Euclidean distance between two independent points
uniformly distributed on $\partial K$ with respect to arclength.
We prove the sharp inequality
\[
    \frac{\theta(K)}{\operatorname{per}K}>\frac{1}{6}.
\]
The constant is approached by rectangles degenerating to a
segment.
\end{abstract}

\end{frontmatter}

\section{Introduction}

Let $K\subset\R^2$ be a convex body, and let $B_1,B_2$ be
independent random points uniformly distributed on $\partial K$
with respect to arclength. Put
\[
    \theta(K)=\E|B_1-B_2|,
\]
and denote the perimeter of $K$ by $\operatorname{per}K$.

According to Tokmachev~\cite[Eq.~(5)]{Tokmachev2026},
Gusakova and Zaporozhets proposed the pair of sharp inequalities
\[
    \frac16<
    \frac{\theta(K)}{\operatorname{per}K}
    \le\frac{2}{\pi^2}.
\]
The upper bound, with equality only for the disk, was proved
by Tokmachev~\cite{TokmachevBoundary}. The purpose of the
present paper is to establish the complementary lower bound.

\begin{thm}\label{thm:lower-boundary-distance}
For every planar convex body $K$,
\[
    \theta(K)>\frac{\operatorname{per}K}{6}.
\]
The constant $1/6$ is sharp.
\end{thm}

Our result is closely related to the Zaporozhets--Tarasov
conjecture, which concerns the comparison of mean distances
between independent uniform points inside a convex body and
on its boundary. Bonnet, Gusakova, Th\"ale and Zaporozhets
\cite{BonnetEtAl} proved that the mean distance $\Delta(K)$
between two independent uniform points in a planar convex body
satisfies
\[
    \Delta(K)<\frac{\operatorname{per}K}{6}.
\]
Consequently, Theorem~\ref{thm:lower-boundary-distance} gives
another proof of the planar Zaporozhets--Tarasov inequality
$\Delta(K)<\theta(K)$. This inequality was previously established
by a different approach in \cite{ShenAug2026}, where a stronger
directional comparison was obtained.

Our proof is based on a comparison of one-dimensional
distributions associated with portions of $\partial K$ visible from a certain direction, which we call \textit{half-boundaries}.
We consider the tangent angle of a uniformly sampled
point on a half-boundary and the median of three independent such angles.
Using the concavity of the Gini mean difference under
mixtures, we first bound $\theta(K)$ from below by an integral
involving the median distribution. A probabilistic comparison
then allows us to replace the median by a single tangent
angle, after which Cauchy's perimeter formula gives the
desired constant $1/6$. Finally, we examine the equality
conditions to establish strictness. 

\par We first recall the central notion and its basic properties.

\begin{defi}
    For a probability measure $\rho$ on $\R$, write
$$
 \Gini(\rho)=\iint_{\R^2}|x-y|\,d\rho(x)d\rho(y).
$$
This is the \textit{Gini mean difference} of $\rho$.
\end{defi}

We use the following standard identities for the Gini mean
difference; see also \cite[Sec.~2.1.1]{Yitzhaki}.

\begin{lem}\label{lem:gini-identities}
Let $\rho,\sigma$ be probability measures on $\R$ with finite
first moments, and let $F,G$ be their distribution functions.
Then
\begin{equation}\label{eq:gini-cdf}
\Gini(\rho)=2\int_\R F(x)(1-F(x))\,dx.
\end{equation}
Moreover, for $0\le\lambda\le1$,
\begin{equation}\label{eq:gini-mixture}
\Gini(\lambda\rho+(1-\lambda)\sigma)
=\lambda\Gini(\rho)+(1-\lambda)\Gini(\sigma)+2\lambda(1-\lambda)\int_\R(F-G)^2\,dx.
\end{equation}
Finally, if $X$ has an atomless distribution $\rho$, then
\begin{equation}\label{eq:gini-cov}
\Gini(\rho)=4\Cov(X,F(X)).
\end{equation}
\end{lem}

\begin{proof}
Let $X,X'$ be independent with law $\rho$. Since
\[
|X-X'|=\int_\R
\left|\mathbf{1}_{\{X\le x\}}-
\mathbf{1}_{\{X'\le x\}}\right|\,dx,
\]
taking expectations and applying Tonelli's theorem gives
\eqref{eq:gini-cdf}.

The distribution function of $\lambda\rho+(1-\lambda)\sigma$
is $H=\lambda F+(1-\lambda)G$. The identity
\[
H(1-H)=\lambda F(1-F)+(1-\lambda)G(1-G)
+\lambda(1-\lambda)(F-G)^2
\]
then proves \eqref{eq:gini-mixture}.

Finally, if $\rho$ is atomless, then $F(X)$ is uniform on
$[0,1]$. By independence and symmetry,
\[
2\Cov(X,F(X))=2\E[X\mathbf{1}_{\{X'\le X\}}]-\E X=\E[(X-X')\mathbf{1}_{\{X'\le X\}}]=\frac12\E|X-X'|.
\]
This proves \eqref{eq:gini-cov}.
\end{proof}

\section{The lower bound}

  The argument we use to prove Theorem~\ref{thm:lower-boundary-distance} is  a
comparison of one-dimensional distributions.  The relevant laws are obtained by
choosing either one random tangent direction or the median of three random tangent
directions on a moving half-boundary.

Orient \(\partial K\) counterclockwise and parametrize it by arclength \(s\in[0,P]\), where \(P=\operatorname{per}K\):
\[
r:[0,P]\longrightarrow\partial K.
\]
\[
r'(s)=(\cos\alpha(s),\sin\alpha(s))
\qquad\text{for almost every }s.
\]
Its tangent angle \(\alpha(s)\) is defined almost everywhere and admits a nondecreasing lift. Let \(\mu\) be the pushforward of arclength under \(s\mapsto\alpha(s)\), and put
\[
P=\mu(S^1),\qquad S(t)=\mu([t,t+\pi]).
\]
For each \(t\), let \(I_t\in[0,\pi]\) have the normalized restriction of \(\mu\) to \([t,t+\pi]\), shifted by \(-t\), and let \(M_t\) be the median of three independent copies of \(I_t\). If \(\partial K\) is \(C^2\) with positive curvature, then \(\mu(d\alpha)=\tau(\alpha)\,d\alpha\), where \(\tau\) is the radius of curvature.

\begin{prop}\label{prop:random-separator}
With the preceding notation,
\begin{equation}\label{eq:median-reduction}
    \theta(K)
    \ge \frac{1}{6P}\int_0^{2\pi}
       S(t)^2\,\mathbb E\sin M_t\,dt.
\end{equation}
\end{prop}

\begin{proof}
Fix $t$ such that $\mu(\{t\})=\mu(\{t+\pi\})=0$.
Such $t$ form a set of full Lebesgue measure, since
$\mu$ has at most countably many atoms. Project the body onto
$$
    u_t=(-\sin t,\cos t).
$$
The boundary $\partial K$ splits into two arcs joining the two points where the projection onto \(u_t\) attains its minimum and maximum. Let $\Gamma_t$ and $\Gamma_{t+\pi}$ be the two boundary arcs whose tangent angles range over $[t, t+\pi]$ and $[t+\pi, t+2\pi]$. Their arclengths are $S(t)$ and $S(t+\pi)$ = $P-S(t)$. Let $\mu_t, \mu_{t+\pi}$ be the pushforwards of normalized arclength measures on these arcs under projection onto $u_t$. Finally, let \(\nu_t\) be the pushforward of normalized arclength measure on the whole boundary \(\partial K\) under projection onto \(u_t\). Since
$$
\partial K=\Gamma_t\cup\Gamma_{t+\pi},
\qquad
S(t)+S(t+\pi)=P,
$$
we have
$$
\nu_t
=
\frac{S(t)}{P}\mu_t
+
\frac{S(t+\pi)}{P}\mu_{t+\pi}.
$$
As the Gini functional is concave under mixtures, 
\begin{equation}\label{concave}
\Gini(\nu_t)\ge\frac{S(t)}{P}\Gini(\mu_t)+\frac{S(t+\pi)}{P}\Gini(\mu_{t+\pi}).\end{equation}
Now
$$
|z|
=
\frac12\int_0^\pi
|\langle z,u_t\rangle|\,dt,
$$
and hence
$$
\theta(K)
=
\frac12\int_0^\pi \Gini(\nu_t)\,dt.
$$
Using \ref{concave},

$$\theta(K) \ge\frac1{2P}\int_0^\pi\left[S(t)\Gini(\mu_t)+S(t+\pi)\Gini(\mu_{t+\pi})\right]dt.$$
Projection of $\Gamma_{t+\pi}$ onto $u_{t+\pi}=-u_t$ reflects the
law obtained by projection onto $u_t$. Since the Gini mean difference
is invariant under reflection, the change of variable $r=t+\pi$ gives
\begin{equation}\label{theta_lower_bound}\theta(K) \geq \frac1{2P}\left[\int_0^\pi S(t)\Gini(\mu_t)\,dt+\int_\pi^{2\pi}S(t)\Gini(\mu_t)\,dt\right]=\frac1{2P}\int_0^{2\pi}S(t)\Gini(\mu_t)\,dt.
\end{equation}

Parametrize $\Gamma_t$ by arclength as
$$
r_t:[0,S(t)]\longrightarrow \Gamma_t,
$$
starting at the point minimizing the projection onto \(u_t\), and put
$$
x_t(s)=\langle r_t(s),u_t\rangle.
$$
If $\alpha_t(s)$ denotes the tangent angle at $r_t(s)$, then
$$
r_t'(s)=(\cos\alpha_t(s),\sin\alpha_t(s))
$$
and therefore
$$
x_t'(s)
=
\langle r_t'(s),u_t\rangle
=
\sin(\alpha_t(s)-t)\ge0.
$$
Thus $x_t$ is nondecreasing on $[0,S(t)]$.  Using that $x_t$ is
nondecreasing and applying Fubini's theorem we get
\begin{equation}\label{eq:chain-energy}
\begin{aligned}
    S(t)^2\Gini(\mu_t)
    &=\int_0^{S(t)}\int_0^{S(t)}
       |x_t(s)-x_t(r)|\,dr\,ds = 2\iint_{0\leq r\leq s \leq S(t)} x_t(s)-x_t(r)\,dr\,ds \\
    &=2\iint_{0\leq r\leq s \leq S(t)} \int_r^s x_t'(z)\,dz\,dr\,ds=2\int_0^{S(t)}z\bigl(S(t)-z\bigr)x_t'(z)\,dz.
\end{aligned}
\end{equation}
Observe that the median $Z_{(2)}$ of three independent uniform points $Z_1, Z_2,Z_3$ on $[0,S(t)]$ has density
$\frac{6z(S(t)-z)}{S(t)^3}\,dz, \,\,\,\text{so}$
    $$S(t)^2\Gini(\mu_t) = \frac{S(t)^3}{3} \E x_t'(Z_{(2)}).$$ Since $s \mapsto \alpha_t(s)$ is nondecreasing on the half-boundary, taking the median commutes with change of variable, i.e.
$$\alpha_t(Z_{(2)}) - t  = M_t$$

and, by definition $x_t'(s) = \sin(\alpha_t(s) - t)$.  Hence
$$    S(t)^2\Gini(\mu_t)=\frac{S(t)^3}{3}\,\mathbb E\sin M_t.$$

Finally, plugging this into \ref{theta_lower_bound} gives
\begin{equation}\label{tehta_lower_bound_fin}
\theta(K)
\ge
\frac1{6P}
\int_0^{2\pi}
S(t)^2\,\E\sin M_t\,dt.
\end{equation}
\end{proof}
Let us replace $M_t$ in the integral above with $I_t$. Then
$$
\begin{aligned}
\int_0^{2\pi}S(t)^2\,\E\sin I_t\,dt
&=
\int_0^{2\pi}
S(t)
\int_t^{t+\pi}
\sin (\omega-t)\,\mu(d\omega)\,dt.
\end{aligned}
$$
The inner integral is the width $w(t)$ of $K$ in the direction $u_t$. Hence
\begin{equation}\label{second_brick}
\int_0^{2\pi}S(t)^2\E\sin I_tdt=\int_0^{2\pi}S(t)w(t)dt=\int_0^\pi\bigl(S(t)+S(t+\pi)\bigr)w(t)dt=P\int_0^\pi w(t)dt=P^2,
\end{equation}
where in the last step we used Cauchy's perimeter formula.
Therefore, the following Proposition together with Proposition~\ref{prop:random-separator} would conclude the proof of Theorem~\ref{thm:lower-boundary-distance} with non-strict inequality.
\begin{prop}\label{prop:median-concentration}
For every $0\le a\le\pi/2$,
\begin{equation}\label{eq:central-concentration}
 \int_0^{2\pi}S(t)^2\mathbb P\{M_t\in[a,\pi-a]\}\,dt
 \ge
 \int_0^{2\pi}S(t)^2\mathbb P\{I_t\in[a,\pi-a]\}\,dt.
\end{equation}
Consequently,
\begin{equation}\label{eq:median-improves-sine}
 \int_0^{2\pi}S(t)^2\mathbb E\sin M_t\,dt
 \ge
 \int_0^{2\pi}S(t)^2\mathbb E\sin I_t\,dt.
\end{equation}
\end{prop}

\begin{proof}
Put $b=\pi-a$ and let, as before, $\mu(d\alpha)$ be the
tangent-angle measure. Set
\[
 C:=\int_0^{2\pi}S(t)^2\,dt,
 \qquad f_T(t):=\frac{S(t)^2}{C},
\]
and choose $T\in[0,2\pi)$ with density $f_T$.
Conditional on $T=t$, choose three independent points
$X,Y,Z$ of the angle circle from normalized $\mu$ on the half-arc
$[t,t+\pi]$. Write $m=\operatorname{med}(X,Y,Z)$, taking real
representatives in this interval. Conditional on $T=t$, the laws of
$X-t$ and $m-t$ are those of $I_t$ and $M_t$, respectively.
Define the events
\[
 E_m:=\{m\in[T+a,T+b]\},
 \qquad E_X:=\{X\in[T+a,T+b]\}.
\]
Total probability, applied to the start $T$, gives
 $$\mathbb P(E_m)=\int_0^{2\pi}\mathbb P(E_m\mid T=t)\,f_T(t)\,dt=\frac1C\int_0^{2\pi}S(t)^2
       \mathbb P\{M_t\in[a,b]\}\,dt,$$
and similarly
\[
 \mathbb P(E_X)
 =\frac1C\int_0^{2\pi}S(t)^2
       \mathbb P\{I_t\in[a,b]\}\,dt.
\]
Thus \eqref{eq:central-concentration} is equivalent to
$\mathbb P(E_m)\ge\mathbb P(E_X)$. We will prove the stronger
conditional inequality
\begin{equation}\label{eq:cond_comparison}
 \mathbb P(E_m\mid X,Y)\ge\mathbb P(E_X\mid X,Y)
 \quad\text{almost surely},
\end{equation}
and then take expectations over the pair.

The joint law of $T,X,Y$ is
$$\frac{S(t)^2}{C}\,dt\,
 \frac{\mathbf 1_{\{x\in[t,t+\pi]\}}\,d\mu(x)}{S(t)}
 \frac{\mathbf 1_{\{y\in[t,t+\pi]\}}\,d\mu(y)}{S(t)}=\frac1C\mathbf 1_{\{x,y\in[t,t+\pi]\}}\,
      dt\,d\mu(x)\,d\mu(y).$$

Hence, conditional on $X,Y$, the start $T$ is uniform among the
half-windows containing both points.

If $X=Y$, then $m=X$, and the conditional comparison ~\eqref{eq:cond_comparison} is an equality.
 Fix a pair $X=x$, $Y=y$, order its
points along their shorter arc, rotate the coordinates and $\mu$,
and lift the points to $0,d$, where $0<d<\pi$. Thus $X,Y$ take the
values $0,d$ in either order. Using the corresponding lift of the
start, conditional on this fixed pair, $T$ is uniform on
\[
 J=[d-\pi,0].
\]
In particular, events with $X$ and $Y$ antipodal occur with probability zero, so we assume $d < \pi$. 

For each $t\in J$, let $Z_t$ have normalized $\mu$-law on $[t,t+\pi]$,
and put
\[
 m_t=\operatorname{med}(0,d,Z_t),\qquad
 p(t)=\mathbb P\{m_t\in[t+a,t+b]\}.
\]
Conditional on $T=t$, the third sample is
independent of $X,Y$, so conditioning on the pair leaves its normalized
$\mu$-law on $[t,t+\pi]$ unchanged. Consequently,
 $$p(t)=\frac1{S(t)}\int_t^{t+\pi}
   \mathbf 1_{\{\operatorname{med}(0,d,z)\in[t+a,t+b]\}}\,\mu(dz)=\mathbb P(E_m\mid X,Y,T=t).$$
 Conditional total probability gives
 $$\mathbb P(E_m\mid X,Y)=\int_J\mathbb P(E_m\mid X,Y,T=t)\,
          \mathbb P\{T\in dt\mid X,Y\}=\frac{1}{\pi-d}\int_Jp(t)\,dt.$$
For either fixed point $q\in\{0,d\}$, the event
$q\in[T+a,T+b]$ is equivalent to $T\in[q-b,q-a]$. Hence
\[
 \begin{aligned}
 \mathbb P\{q\in[T+a,T+b]\mid X,Y\}
 &=\frac1{\pi-d}\int_J
     \mathbf 1_{\{q\in[t+a,t+b]\}}\,dt\\
 &=\frac{|J\cap[q-b,q-a]|}{\pi-d}
 =\frac{\delta}{\pi-d},
 \end{aligned}
\]
where the two interval intersections have the same length
$\delta=\max\{0,\min\{b-a,b-d\}\}$.
Since $X$ is one of these two fixed points,
$\mathbb P(E_X\mid X,Y)=\delta/(\pi-d)$.
Thus the conditional comparison reduces to proving
\[
 \int_J p(t)\,dt\ge\delta.
\]
If $\delta=0$, this follows from $p(t)\ge0$. Assume $\delta>0$.

We estimate the probability that $m_t$ falls outside $[t+a,t+b]$,
on either its left or its right. The key observation is that the law
of $m_t$ increases stochastically with $t$. Indeed, for $0<c\le d$,
\[
 \mathbb P\{m_t<c\}
 =\frac{\mu[t,c)}{\mu[t,c)+\mu[c,t+\pi]}.
\]
As $t$ increases, $\mu[t,c)$ decreases and $\mu[c,t+\pi]$ increases,
so the displayed probability is nonincreasing. For $c\le0$ or $c>d$,
the probability is constant.

Write
\[
 A(t)=\mathbb P\{m_t<t+a\},\qquad
 B(t)=\mathbb P\{m_t>t+b\}.
\]
These are the probabilities of the two disjoint ways to miss the
interval, so
\[
 p(t)=1-A(t)-B(t).
\]
Since $0\le m_t\le d$ and $\delta\le b-d$,
\[
 \begin{aligned}
 A(t)&=0 &&(d-\pi\le t\le d-\pi+\delta),\\
 B(t)&=0 &&(-\delta\le t\le0).
 \end{aligned}
\]
Moreover, $\delta\le b-a$ and stochastic monotonicity give, for
$d-\pi\le t\le-\delta$,
\[
 A(t+\delta)
 \le\mathbb P\{m_t<t+\delta+a\}
 \le\mathbb P\{m_t\le t+b\}
 =1-B(t).
\]
To integrate this comparison, first remove the intervals on which
$A$ and $B$ vanish. For $A$, this gives
$$\int_{d-\pi}^{0}A(u)\,du=\int_{d-\pi+\delta}^{0}A(u)\,du=\int_{d-\pi}^{-\delta}A(t+\delta)\,dt,$$
 For $B$, we have
\[
 \int_{d-\pi}^{0}B(t)\,dt
 =\int_{d-\pi}^{-\delta}B(t)\,dt.
\]
 Summing and using $A(t+\delta)+B(t)\le1$, we obtain
 $$\int_J\bigl(A(t)+B(t)\bigr)\,dt=\int_{d-\pi}^{-\delta}
      \bigl(A(t+\delta)+B(t)\bigr)\,dt\le\int_{d-\pi}^{-\delta}1\,dt
 =\pi-d-\delta.$$
Since $p(t)=1-A(t)-B(t)$, it follows that
 $$\int_Jp(t)\,dt=|J|-\int_J\bigl(A(t)+B(t)\bigr)\,dt\ge(\pi-d)-(\pi-d-\delta)=\delta.$$
Dividing by $\pi-d$ proves the conditional comparison:
\[
 \mathbb P(E_m\mid X,Y)=\frac1{\pi-d}\int_Jp(t)\,dt\ge\frac{\delta}{\pi-d}
 =\mathbb P(E_X\mid X,Y).
\]
Finally, taking expectations with respect to the joint law of $X,Y$
and using total probability yields
$$\int_0^{2\pi}S(t)^2\mathbb P\{M_t\in[a,b]\}\,dt=C\,\mathbb P(E_m)=C\,\mathbb E\bigl[\mathbb P(E_m\mid X,Y)\bigr]$$ $$ \ge C\,\mathbb E\bigl[\mathbb P(E_X\mid X,Y)\bigr]=C\,\mathbb P(E_X)=\int_0^{2\pi}S(t)^2\mathbb P\{I_t\in[a,b]\}\,dt.$$
As $b=\pi-a$, this is \eqref{eq:central-concentration}.
This argument also applies to non-smooth
tangent-angle measures; choices at window endpoints affect only a
Lebesgue-null set of starts.

The identity
\[
 \sin r=\int_0^{\pi/2}\mathbf 1_{[a,\pi-a]}(r)\cos a\,da,
 \qquad 0\le r\le\pi,
\]
and Tonelli's theorem yield \eqref{eq:median-improves-sine} by integrating
\eqref{eq:central-concentration} against $\cos a\,da$.
\end{proof}

\begin{proof}[Proof of Theorem~\ref{thm:lower-boundary-distance}]
As we outlined before Proposition~\ref{prop:median-concentration},
\begin{equation}\label{eq:one-sample-baseline}
\begin{aligned}
    \int_0^{2\pi}S(t)^2\mathbb E\sin I_t\,dt
    &=\int_0^{2\pi}S(t)
      \left(\int_t^{t+\pi}\sin(\alpha-t)\,\mu(d\alpha)\right)dt\\
    &=\int_0^{2\pi}S(t)w(t)\,dt,
\end{aligned}
\end{equation}
where $w(t)$ is the width of $K$ in direction $u_t$. Indeed, $$w(t)=\langle r_t(S(t))-r_t(0),u_t\rangle=\int_0^{S(t)}\langle r_t'(s),u_t\rangle\,ds=\int_0^{S(t)}\sin(\alpha_t(s)-t)\,ds$$$$=\int_t^{t+\pi}\sin(\alpha-t)\,\mu(d\alpha).$$ Since
$$
    S(t)+S(t+\pi)=P,
    \qquad
    w(t+\pi)=w(t),
$$

Cauchy's perimeter formula gives
\begin{equation}\label{eq:baseline-p-squared}
    \int_0^{2\pi}S(t)w(t)\,dt
    =P\int_0^\pi w(t)\,dt=P^2.
\end{equation}
Combining Propositions~\ref{prop:random-separator} and ~\ref{prop:median-concentration} with
\eqref{eq:baseline-p-squared}, we obtain
$$
    \theta(K)\ge \frac{P}{6}.
$$
 It remains to exclude equality for a nondegenerate convex body.
Suppose therefore that
\[
    \theta(K)=\frac{P}{6}.
\]
Then both inequalities supplied by Propositions~\ref{prop:median-concentration} and~\ref{prop:random-separator} must be
equalities.

Recall that we denote $\mu$ the tangent-angle measure of $\partial K$, of total mass
$P$. We show that equality in \eqref{eq:median-improves-sine} requires
that no open semicircle contain three points of $\operatorname{supp}\mu$.
Suppose otherwise. Choose small arcs $U,V,W$ about three such points,
lifted in this order to $\mathbb R$, and $a_0>0$ so that
\[
 v-x>a_0,\qquad y-v>a_0,\qquad y-x<\pi-a_0
 \quad(x\in U,\ v\in V,\ y\in W).
\]
All three arcs have positive $\mu$-measure.

Fix $0<a<a_0$ and a pair $x\in U$, $y\in W$. Translate the angular
coordinates so that $x=0$, $y=d$, and use $m_t,p(t),A(t),B(t)$ from the
proof of Proposition~\ref{prop:median-concentration}. Here
$\delta=\pi-d-a>0$. The difference between the two sides of its
integral comparison is
\[
 \begin{aligned}
 \int_{d-\pi}^{0}p(t)\,dt-\delta
 &=\int_{d-\pi}^{d-\pi+a}
   \bigl[1-B(t)-A(t+\delta)\bigr]\,dt.
 \end{aligned}
\]
For $t=d-\pi+s$ with $0<s<a$, stochastic monotonicity gives
\[
 \begin{aligned}
 1-B(t)-A(t+\delta)
 &=\mathbb P\{m_t\le d+s-a\}
   -\mathbb P\{m_{t+\delta}<s\}\\
 &\ge\mathbb P\{s\le m_t\le d+s-a\}
 \ge\frac{\mu(V)}{P}>0.
 \end{aligned}
\]
Indeed, the translated arc $V$ lies in $(a,d-a)$, hence in
$[s,d+s-a]$, and a third sample in $V$ is the median. The excess is
therefore at least $a\mu(V)/P$ for every pair in $U\times W$.
The pair disintegration in Proposition~\ref{prop:median-concentration}
now makes \eqref{eq:central-concentration} strict for every $0<a<a_0$.
Integration against $\cos a\,da$ makes
\eqref{eq:median-improves-sine} strict, a contradiction.
Now, we have shown that
\[
    |\operatorname{supp}\mu|\le4.
\]
Indeed, if five points occurred in cyclic order and $g_1,\ldots,g_5$
were the intervening angular gaps, then the absence of three points
in an open semicircle would imply
\[
    g_i+g_{i+1}\ge\pi
    \qquad (i\!\!\!\pmod 5).
\]
Summing gives $4\pi\ge5\pi$, since $\sum_i g_i=2\pi$, which is
impossible.

Hence $K$ is a polygon with at most four sides. We now show that the inequality of 
Proposition~2.1 is strict for every nondegenerate polygon. Choose an
open interval $J$ of angles $t$ so that for every $t\in J$
$$u_t = (-\sin t, \cos t)$$
is not normal to any edge of $K$.
For $t\in J$ the linear functional $x\mapsto\langle x, u_t \rangle$ attains its minimum and maximum at unique vertices $v$ and $w$, respectively. Those vertices split $\partial K$ in two complementary polygonal chains. Denote their lengths by $S_1,S_2$, and let $\nu_{1,t},\nu_{2,t}$ be their normalized arclength measures projected onto $u_t$.

The equality in the concavity step of
Proposition~2.1 is equivalent to
\[
    \nu_{1,t}=\nu_{2,t}.
\]

Let $e_1,e_2$ be the unit vectors along the edges of
the two chains incident to $v$ and directed away from $v$. For $t\in J$,
\[
    \langle e_i,u_t\rangle>0.
\]
Near the projection of $v$, the densities of
$\nu_{1,t}$ and $\nu_{2,t}$ are respectively
\[
    \frac{1}{S_1\langle e_1,u_t\rangle},
    \qquad
    \frac{1}{S_2\langle e_2,u_t\rangle}.
\]
Therefore equality of the two projected laws would imply
\[
    \bigl\langle S_1e_1-S_2e_2,u_t\bigr\rangle=0.
\]
As $S_1e_1-S_2e_2\ne0$, this can hold only for isolated values of $t$, whereas $J$ is an interval. Hence the Gini
concavity inequality is strict on a set of positive measure, and so
Proposition~2.1 is strict for every nondegenerate polygon.

Therefore,
\[
    \theta(K)>\frac{\operatorname{per}K}{6}.
\]
Sharpness follows from the rectangles
$K_\varepsilon=[0,1]\times[0,\varepsilon]$.
As $\varepsilon\downarrow0$, normalized boundary arclength concentrates
on the two horizontal sides, and hence
$\theta(K_\varepsilon)\to\mathbb E|U-V|=1/3$, where
$U,V\sim\operatorname{Unif}[0,1]$, while
$\operatorname{per}K_\varepsilon\to2$.
\end{proof}

\end{document}